# Statistical Convergence and Convergence in Statistics

## Mark Burgin[a], Oktay Duman[b]

[a]*Department of Mathematics, University of California, Los Angeles, California 90095-1555, USA*

[b]*TOBB Economics and Technology University, Faculty of Arts and Sciences, Department of Mathematics, Sögütözü 06530, Ankara, Turkey*

## Abstract

Statistical convergence was introduced in connection with problems of series summation. The main idea of the statistical convergence of a sequence $l$ is that the majority of elements from $l$ converge and we do not care what is going on with other elements. We show (Section 2) that being mathematically formalized the concept of statistical convergence is directly connected to convergence of such statistical characteristics as the mean and standard deviation. At the same time, it known that sequences that come from real life sources, such as measurement and computation, do not allow, in a general case, to test whether they converge or statistically converge in the strict mathematical sense. To overcome limitations induced by vagueness and uncertainty of real life data, neoclassical analysis has been developed. It extends the scope and results of the classical mathematical analysis by applying fuzzy logic to conventional mathematical objects, such as functions, sequences, and series. The goal of this work is the further development of neoclassical analysis. This allows us to reflect and model vagueness and uncertainty of our knowledge, which results from imprecision of measurement and inaccuracy of computation. In the context on the theory of fuzzy limits, we develop the structure of statistical fuzzy convergence and study its properties. Relations between statistical fuzzy convergence and fuzzy convergence are considered in Theorems 3.1 and 3.2. Algebraic structures of statistical fuzzy limits are described in Theorem 3.5. Topological structures of statistical fuzzy limits are described in Theorems 3.3 and 3.4. Relations between statistical fuzzy convergence and fuzzy convergence of statistical characteristics, such as the mean (average) and standard deviation, are studied in Section 4. Introduced constructions and obtained results open new directions for further research that are considered in the Conclusion.

## 1. Introduction

The idea of statistical convergence goes back to the first edition (published in Warsaw in 1935) of the monograph of Zygmund [37]. Formally the concept of statistical convergence was introduced by Steinhaus [34] and Fast [18] and later reintroduced by Schoenberg [33].

Statistical convergence, while introduced over nearly fifty years ago, has only recently become an area of active research. Different mathematicians studied properties of statistical convergence and applied this concept in various areas such as measure theory [30], trigonometric series [37], approximation theory [16], locally convex spaces [29], finitely additive set functions [14], in the study of subsets of the Stone-Chech compactification of the set of natural numbers [13], and Banach spaces [15].

However, in a general case, neither limits nor statistical limits can be calculated or measured with absolute precision. To reflect this imprecision and to model it by mathematical structures, several approaches in mathematics have been developed: fuzzy set theory, fuzzy logic, interval analysis, set valued analysis, etc. One of these approaches is the neoclassical analysis (cf., for example, [7, 8]). In it, ordinary structures of analysis, that is, functions, sequences, series, and operators, are studied by means of fuzzy concepts: fuzzy limits, fuzzy continuity, and fuzzy derivatives. For example, continuous functions, which are studied in the classical analysis, become a part of the set of the fuzzy continuous functions studied in neoclassical analysis. Neoclassical analysis extends methods of classical calculus to reflect uncertainties that arise in computations and measurements.

The aim of the present paper is to extend and study the concept of statistical convergence utilizing a fuzzy logic approach and principles of the *neoclassical analysis*, which is a new branch of fuzzy mathematics and extends possibilities provided by the classical analysis [7, 8]. Ideas of fuzzy logic have been used not only in many applications, such as, in bifurcation of non-linear dynamical systems, in the control of chaos, in the computer programming, in the quantum physics, but also in various branches of mathematics, such as, theory of metric and topological spaces, studies of convergence of sequences and functions, in the theory of linear systems, etc.

In the second section of this paper, going after introduction, we remind basic constructions from the theory of statistical convergence consider relations between statistical convergence, ergodic systems, and convergence of statistical characteristics such as the mean (average), and standard deviation. In the third section, we introduce a new type of fuzzy convergence, the concept of statistical fuzzy convergence, and give a useful characterization of this type of convergence. In the fourth section, we consider relations between statistical fuzzy convergence and fuzzy convergence of statistical characteristics such as the mean (average) and standard deviation.

For simplicity, we consider here only sequences of real numbers. However, in a similar way, it is possible to define statistical fuzzy convergence for sequences of complex numbers and obtain similar properties.

## 2. Convergence in statistics

Statistics is concerned with the collection and analysis of data and with making estimations and predictions from the data. Typically two branches of statistics are discerned: descriptive and inferential. Inferential statistics is usually used for two tasks: to estimate properties of a population given sample characteristics and to predict properties of a system given its past and current properties. To do this, specific statistical constructions were invented. The most popular and useful of them are the average or mean (or more exactly, arithmetic mean) $\mu$ and standard deviation $\sigma$ (variance $\sigma^2$).

To make predictions for future, statistics accumulates data for some period of time. To know about the whole population, samples are used. Normally such inferences (for future or for population) are based on some assumptions on limit processes and their convergence. Iterative processes are used widely in statistics. For instance the empirical approach to probability is based on the law (or better to say, conjecture) of big numbers, states that a procedure repeated again and again, the relative frequency probability tends to approach the actual probability. The foundation for estimating population parameters and hypothesis testing is formed by the central limit theorem, which tells us how sample means change when the sample size grows. In experiments, scientists measure how statistical characteristics (e.g., means or standard deviations) converge (cf., for example, [23, 31]).

Convergence of means/averages and standard deviations have been studied by many authors and applied to different problems (cf. [1-4, 17, 19, 20, 24-28, 35]). Convergence of statistical characteristics such as the average/mean and standard deviation are related to statistical convergence as we show in this section and Section 4.

Consider a subset $K$ of the set $N$ of all natural numbers. Then $K_n = \{k \in K; k \leq n\}$.

**Definition 2.1.** The asymptotic density $d(K)$ of the set $K$ is equal to

$$\lim_{n \to \infty} (1/n) \, |K_n|$$

whenever the limit exists; here $|B|$ denotes the cardinality of the set $B$.

Let us consider a sequence $l = \{a_i ; i = 1,2,3,\ldots\}$ of real numbers, real number $a$, and the set $L_\varepsilon(a) = \{i \in N; |a_i - a| \geq \varepsilon\}$.

**Definition 2.2.** The asymptotic density, or simply, density $d(l)$ of the sequence $l$ with respect to $a$ and $\varepsilon$ is equal to $d(L_e(a))$.

Asymptotic density allows us to define statistical convergence.

**Definition 2.3.** A sequence $l = \{a_i ; i = 1, 2, 3, \ldots\}$ is statistically convergent to $a$ if $d(L_e(a)) = 0$ for every $\varepsilon > 0$. The number (point $a$) is called the statistical limit of $l$. We denote this by $a = stat\text{-}\lim l$.

Note that convergent sequences are statistically convergent since all finite subsets of the natural numbers have density zero. However, its converse is not true [21, 33]. This is also demonstrated by the following example.

**Example 2.1.** Let us consider the sequence $l = \{a_i \, ; \, i = 1,2,3,\dots\}$ whose terms are

$$a_i = \begin{cases} i & \text{when } \; i = n^2 \text{ for all } n = 1,2,3,\dots \\ \\ 1/i & \text{otherwise.} \end{cases}$$

Then, it is easy to see that the sequence $l$ is divergent in the ordinary sense, while $0$ is the statistical limit of $l$ since $d(K) = 0$ where $K = \{n^2 \text{ for all } n = 1,2,3,\dots\}$.

Not all properties of convergent sequences are true for statistical convergence. For instance, it is known that a subsequence of a convergent sequence is convergent. However, for statistical convergence this is not true. Indeed, the sequence $h = \{i \, ; \, i = 1,2,3,\dots\}$ is a subsequence of the statistically convergent sequence $l$ from Example 2.1. At the same time, $h$ is statistically divergent.

However, if we consider dense subsequences of fuzzy convergent sequences, it is possible to prove the corresponding result.

**Definition 2.4.** A subset $K$ of the set $N$ is called *statistically dense* if $d(K) = 1$.

**Example 2.2.** The set $\{ \, i \neq n^2 \, ; \, i = 1,2,3,\dots; \, n = 1,2,3,\dots \}$ is statistically dense, while the set $\{ \, 3i; \, i = 1,2,3,\dots \}$ is not.

**Lemma 2.1. a)** A statistically dense subset of a statistically dense set is a statistically dense set.

**b)** The intersection and union of two statistically dense sets are statistically dense sets.

**Definition 2.5.** A subsequence $h$ of the sequence $l$ is called *statistically dense* in $l$ if the set of all indices of elements from $h$ is statistically dense.

**Corollary 2.1. a)** A statistically dense subsequence of a statistically dense subsequence of *l* is a statistically dense subsequence of *l*.

**b)** The intersection and union of two statistically dense subsequences are statistically dense subsequences.

**Theorem 2.1.** A sequence *l* is statistically convergent if and only if any statistically dense subsequence of *l* is statistically convergent.

*Proof*. *Necessity*. Let us take a statistically convergent sequence $l = \{a_i\,;\, i = 1,2,3,\ldots\}$ and a statistically dense subsequence $h = \{b_k\,;\, k = 1,2,3,\ldots\}$ of *l*. Let us also assume that *h* statistically diverges. Then for any real number *a*, there is some $\varepsilon > 0$ such that $liminf_{n \to \infty} (1/n)\,|H_{n,\varepsilon}(a)| = d > 0$ for some $d \in (0, 1)$, where $H_{n,\varepsilon}(a) = \{k \le n;\, |b_k - a| > \varepsilon\}$. As *h* is a subsequence of *l*, we have $L_{n,\varepsilon}(a) \supseteq H_{n,\varepsilon}(a)$ where $L_{n,\varepsilon}(a) = \{i \le n;\, |a_i - a| > \varepsilon\}$. Consequently, $liminf_{n \to \infty} (1/n)|L_{n,\varepsilon}(a)| \ge d > 0$, which yields that $d(\{i;\, |a_i - a| > \varepsilon\}) \ne 0$. Thus, *l* is also statistically divergent.

*Sufficiency* follows from the fact that *l* is a statistically dense subsequence of itself.

**Corollary 2.1.** A statistically dense subsequence of a statistically convergent sequence is statistically convergent.

To each sequence $l = \{a_i\,;\, i = 1,2,3,\ldots\}$ of real numbers, it is possible to correspond a new sequence $\mu(l) = \{\mu_n = (1/n)\, \Sigma_{i=1}^{n}\, a_i\,;\, n = 1,2,3,\ldots\}$ of its partial averages (means). Here a partial average of *l* is equal to $\mu_n = (1/n)\, \Sigma_{i=1}^{n}\, a_i$ .

Sequences of partial averages/means play an important role in the theory of ergodic systems [5]. Indeed, the definition of an ergodic system is based on the concept of the "time average" of the values of some appropriate function *g* arguments for which are dynamic transformations *T* of a point *x* from the manifold of the dynamical system. This average is given by the formula

$$g(x) = \lim (1/n)\, \Sigma_{k=1}^{n-1}\, g(T^k x).$$

In other words, the dynamic average is the limit of the partial averages/means of the sequence $\{ T^k x \, ; k = 1,2,3,\ldots \}$.

Let $l = \{ a_i \, ; \, i = 1,2,3,\ldots \}$ be a bounded sequence, i.e., there is a number $m$ such that $|a_i| < m$ for all $i \in \mathbf{N}$. This condition is usually true for all sequences generated by measurements or computations, i.e., for all sequences of data that come from real life.

**Theorem 2.2.** If $a = stat\text{-}\lim l$, then $a = \lim \mu(l)$.

*Proof.* Since $a = stat\text{-}\lim l$, for every $\varepsilon > 0$, we have

(2.1) $\qquad\qquad \lim_{n \to \infty} (1/n) \, |\{ i \leq n, \, i \in \mathbf{N}; \, |a_i - a| \geq \varepsilon \}| = 0.$

As $|a_i| < m$ for all $i \in \mathbf{N}$, there is a number $k$ such that $|a_i - a| < k$ for all $i \in \mathbf{N}$. Namely, $| a_i - a| \leq | a_i | + | a| \leq m + | a| = k$. Taking the set $L_{n,\varepsilon} \, (a) = \{ i \leq n, \, i \in \mathbf{N}; \, | a_i - a| \geq \varepsilon \}$, denoting $|L_{n,\varepsilon} \, (a)|$ by $u_n$, and using the hypothesis $|a_i| < m$ for all $i \in \mathbf{N}$, we have the following system of inequalities:

$$|\mu_n - a| = |(1/n) \, \Sigma_{i=1}^{n} \, a_i - a|$$

$$\leq (1/n) \, \Sigma_{i=1}^{n} \, | a_i - a|$$

$$\leq (1/n) \, \{ k u_n + ( n - u_n) \varepsilon \}$$

$$\leq (1/n) \, (k u_n + n \varepsilon )$$

$$= \varepsilon + k \, (u_n/n).$$

From the equality (2.1), we get, for sufficiently large $n$, the inequality $|\mu_n - a| < \varepsilon + k \, \varepsilon$. Thus, $a = \lim \mu(l)$.

Theorem is proved.

**Remark 2.1.** However, convergence of the partial averages/means of a sequence does not imply statistical convergence of this sequence as the following example demonstrates.

**Example 2.3.** Let us consider the sequence $l = \{ a_i \, ; \, i = 1,2,3,\ldots \}$ whose terms are $a_i = (-1)^i \sqrt{i} \,$. This sequence is statistically divergent although $\lim \mu(l) = 0$.

Taking a sequence $l = \{a_i ; i = 1,2,3,\dots\}$ of real numbers, it is possible to construct not only the sequence $\mu(l) = \{\mu_n = (1/n) \Sigma_{i=1}^{n} a_i ; n = 1,2,3,\dots\}$ of its partial averages (means) but also the sequences $\sigma(l) = \{\sigma_n = ((1/n) \Sigma_{i=1}^{n} (a_i - \mu_n)^2)^{1/2} ; n = 1,2,3,\dots\}$ of its partial standard deviations $\sigma_n$ and $\sigma^2(l) = \{\sigma_n^2 = (1/n) \Sigma_{i=1}^{n} (a_i - \mu_n)^2 ; n = 1,2,3,\dots\}$ of its partial variances $\sigma_n^2$.

**Theorem 2.3.** If $a = stat\text{-}\lim l$ and $| a_i | < m$ for all $i \in N$, then $\lim \sigma(l) = 0$.

*Proof.* We will show that $\lim \sigma^2(l) = 0$. By the definition, $\sigma_n^2 = (1/n) \Sigma_{i=1}^{n} (a_i - \mu_n)^2 = (1/n) \Sigma_{i=1}^{n} (a_i)^2 - \mu_n^2$. Thus, $\lim \sigma^2(l) = \lim_{n\to\infty} (1/n) \Sigma_{i=1}^{n} (a_i)^2 - \lim_{n\to\infty} \mu_n^2$.

If $| a_i | < m$ for all $i \in N$, then there is a number $p$ such that $| a_i^2 - a^2 | < p$ for all $i \in N$. Namely, $| a_i^2 - a^2 | \leq | a_i |^2 + | a |^2 < m^2 + | a |^2 < m^2 + | a |^2 + m + | a | = p$. Let us consider the absolute value of the difference $\mu_n^2 - (1/n) \Sigma_{i=1}^{n} (a_i)^2 = \sigma_n^2$. Taking the set $L_{n,\varepsilon}(a) = \{i \leq n, i \in N; | a_i - a | \geq \varepsilon \}$, denoting $|L_{n,\varepsilon}(a)|$ by $u_n$, and using the hypothesis $|a_i| < m$ for all $i \in N$, we have the following system of inequalities:

$$|\sigma^2_n | = | (1/n) \Sigma_{i=1}^{n} (a_i)^2 - \mu_n^2 |$$

$$= | (1/n) \Sigma_{i=1}^{n} (a_i^2 - a^2) - (\mu_n^2 - a^2)|$$

$$\leq (1/n) \Sigma_{i=1}^{n} | a_i^2 - a^2 | + |\mu_n^2 - a^2|$$

$$< (p/n) \Sigma_{i=1}^{n} |a_i - a| + |\mu_n^2 - a^2|$$

$$< (p/n) (u_n + (n - u_n)\varepsilon) + |\mu_n^2 - a^2|$$

$$< (p/n) (u_n + n\varepsilon) + |\mu_n^2 - a^2|$$

$$= p (u_n /n) + \varepsilon p + |\mu_n^2 - a^2|$$

as $|a_i^2 - a^2| \leq |a_i - a| \ |a_{i+} a| < |a_i - a| \cdot p$. By Theorem 2.2, we have $a = \lim \mu(l)$, which guarantees that $\lim \mu_n^2 = a^2$. Also by (2.1) $\lim (u_n /n) = 0$. Since $\varepsilon > 0$ was arbitrary, the right hand side of the above inequality tends to zero as $n \to \infty$. Therefore, we have $\lim \sigma(l) = 0$.

Theorem is proved.

**Corollary 2.2.** If $a = stat\text{-}\lim l$ and $|a_i| < m$ for all $i \in \mathbf{N}$, then $\lim \sigma^2(l) = 0$.

**Theorem 2.4.** A sequence $l$ is statistically convergent if its sequence of partial averages $\mu(l)$ converges and $a_i \le \lim \mu(l)$ (or $a_i \ge \lim \mu(l)$) for all $i = 1, 2, 3, \ldots$ .

*Proof*. Let us assume that $a = \lim \mu(l)$, $a_i \le \lim \mu(l)$ and take some $\varepsilon > 0$, the set $L_{n,\varepsilon}(a) = \{i \le n, i \in \mathbf{N}; |a_i - a| \ge \varepsilon\}$, and denote $|L_{n,\varepsilon}(a)|$ by $u_n$ . Then we have

$$|a - \mu_n| = |a - (1/n) \Sigma_{i=1}^{n} a_i|$$

$$= |(1/n) \Sigma_{i=1}^{n} (a - a_i)|$$

$$= (1/n) \Sigma_{i=1}^{n} (a - a_i)$$

$$\ge (1/n) \Sigma_{|a_i - a| \ge \varepsilon} (a - a_i)$$

$$\ge (u_n/n)\varepsilon$$

Consequently, $\lim_{n \to \infty} |a - \mu_n| \ge \lim_{n \to \infty} (u_n/n) \varepsilon$. As $\lim_{n \to \infty} |a - \mu_n| = 0$ and $\varepsilon$ is a fixed number, we have $\lim_{n \to \infty} (1/n) |\{i \le n, i \in \mathbf{N}; |a_i - a| \ge \varepsilon\}| = 0$, i.e., $a = stat\text{-}\lim l$.

The case when $a_i \ge \lim \mu(l)$ for all $i = 1, 2, 3, \ldots$ is considered in a similar way.

Theorem is proved as $\varepsilon$ is an arbitrary positive number.

Let $l = \{a_i ; i = 1,2,3,\ldots\}$ be a bounded sequence, i.e., there is a number $m$ such that $|a_i| < m$ for all $i \in \mathbf{N}$.

**Theorem 2.5.** A sequence $l$ is statistically convergent if and only if its sequence of partial averages $\mu(l)$ converges and its sequence of partial standard deviations $\sigma(l)$ converges to 0.

*Proof*. *Necessity* follows from Theorems 2.2 and 2.3.

*Sufficiency.* Let us assume that $a = \lim \mu(l)$, $\lim \sigma(l) = 0$, and take some $\varepsilon > 0$. This implies that for any $\lambda > 0$, there is a number $n$ such that $\lambda > \mid a - \mu_n \mid$. Then taking a number $n$ such that it implies the inequality $\varepsilon > \lambda$, we have

$$\sigma_n^2 = (1/n) \Sigma_{i=1}^{n} (a_i - \mu_n)^2$$

$$\geq (1/n) \sum \{(a_i - \mu_n)^2 ; |a_i - a| \geq \varepsilon \}$$

$$= (1/n) \sum \{((a_i - a) + (a - \mu_n))^2 ; |a_i - a| \geq \varepsilon \}$$

(2.2) $$> (1/n) \sum \{((a_i - a) \pm \lambda)^2 ; |a_i - a| \geq \varepsilon \}$$

$$= (1/n) \sum \{((a_i - a)^2 \pm 2\lambda(a_i - a) + \lambda^2); |a_i - a| \geq \varepsilon \}$$

(2.3) $$= (1/n) \sum \{(a_i - \mu_n)^2; |a_i - a| \geq \varepsilon \} \pm 2\lambda (1/n) \sum \{(a_i - a); |a_i - a| \geq \varepsilon \} + \lambda^2$$

as $(a_i - \mu_n) = (a_i - a) + (a - \mu_n)$ and we take $+ \lambda$ or $- \lambda$ in the expression (2.2) according to the following rules:

1) if $(a_i - a) \geq 0$ and $(a - \mu_n) \geq 0$, then $(a_i - a) + (a - \mu_n) \geq (a_i - a) > (a_i - a) - \lambda$, and we take $- \lambda$;

2) if $(a_i - a) \geq 0$ and $(a - \mu_n) \leq 0$, then $(a_i - a) + (a - \mu_n) \geq (a_i - a) - |a - \mu_n| > (a_i - a) - \lambda$, and we take $- \lambda$;

3) if $(a_i - a) \leq 0$ and $(a - \mu_n) \geq 0$, then $|(a_i - a) + (a - \mu_n)| = |(a - a_i) - (a - \mu_n)| > |(a - a_i) - \lambda| = |(a_i - a) + \lambda|$, and we take $+ \lambda$;

4) if $(a_i - a) \leq 0$ and $(a - \mu_n) \leq 0$, then $|(a_i - a) + (a - \mu_n)| \geq \mid a - a_i \mid > |(a_i - a) + \lambda|$ as $a_i - a < - \varepsilon$, and we take $+ \lambda$.

In the expression (2.3), it is possible to take a sequence $\{\lambda_k ; k = 1, 2, 3, \dots \}$ such that the sequence $\lambda_k^2$ converges to 0 because the sequence $\{\mu_n ; n = 1,2,3,\dots \}$ converges to $a$ when $n$ tends to $\infty$. The sum $2\lambda_k (1/n) \sum \{(a_i - a); |a_i - a| \geq \varepsilon \}$ also converges to 0 when $k$ tends to $\infty$ because $\lambda_k$ converges to 0 and $(1/n) \sum (a_i - a) < (1/n) \sum (|a_i| + |a|) \leq m + |a|$. At the same time, the sequence $\{\sigma_n; n = 1, 2, 3,\dots \}$ also converges to 0. Thus, $\lim_{n \to \infty}$

$(1/n) \sum \{(a_i - \mu_n)^2; |a_i - a| \geq \varepsilon\} = 0$. This implies that $\lim_{n\to\infty} (1/n) \sum \{|a_i - \mu_n|^2; |a_i - a| \geq \varepsilon\} = 0$. At the same time, $\lim_{n\to\infty} (1/n) \sum \{|a_i - \mu_n|^2; |a_i - a| \geq \varepsilon\} \geq \varepsilon \cdot (\lim_{n\to\infty} (1/n) |\{i \leq n, i \in \mathbf{N}; |a_i - a| \geq \varepsilon\}|)$. As $\varepsilon$ is a fixed number, we have $\lim_{n\to\infty} (1/n) |\{i \leq n, i \in \mathbf{N}; |a_i - a| \geq \varepsilon\}| = 0$ for any $\varepsilon > 0$ as $\varepsilon$ is an arbitrary positive number, i.e., $a = stat\text{-}\lim l$.

Theorem is proved.

### 3. Statistical fuzzy convergence

Here we extend statistical convergence to statistical fuzzy convergence, which, as we have discussed, is more realistic for real life applications.

For convenience, throughout the paper, $r$ denotes a non-negative real number and $l = \{a_i; i = 1,2,3,\ldots\}$ represents a sequence of real numbers.

Let us consider the set $L_{r,\varepsilon}(a) = \{i \in \mathbf{N}; |a_i - a| \geq r + \varepsilon\}$ and a non-negative real number $r \geq 0$.

**Definition 3.1.** A sequence $l$ is statistically $r$-converges to a number a if $d(L_{r,\varepsilon}(a)) = 0$ for every $\varepsilon > 0$. The number (point $a$) is called a statistical $r$-limit of $l$. We denote this by $a = stat\text{-}r\text{-}\lim l$.

Definition 3.1 implies the following results.

**Lemma 3.1. a)** $a = stat\text{-}r\text{-}\lim l \Leftrightarrow \forall \varepsilon > 0$, $\lim_{n\to\infty} (1/n) |\{i \in \mathbf{N}; i \leq n ; |a_i - a| \geq r + \varepsilon\}| = 0$.

**b)** $a = stat\text{-}r\text{-}\lim l \Leftrightarrow \forall \varepsilon > 0$, $\lim_{n\to\infty} (1/n) |\{i \in \mathbf{N}; i \leq n ; |a_i - a| < r + \varepsilon\}| = 1$.

**Remark 3.1.** We know from [8] that if $a = \lim l$ (in the ordinary sense), then for any $r \geq 0$, we have $a = r\text{-}\lim l$. In a similar way, using Definition 3.1, we easily see that if $a = r\text{-}\lim l$, then we have $a = stat\text{-}r\text{-}\lim l$ since every finite subset of the natural numbers has density zero. However, its converse is not true as the following example of a sequence that is statistically $r$-convergent but not $r$-convergent and also not statistically convergent. shows.

**Example 3.1.** Let us consider the sequence $l = \{a_i \, ; i = 1,2,3,\dots\}$ whose terms are

$$(3.1) \qquad a_i = \begin{cases} i & \text{when } i = n^2 \text{ for all } n = 1,2,3,\dots \\ \\ (\text{-}1)^i & \text{otherwise.} \end{cases}$$

Then, it is easy to see that the sequence $l$ is divergent in the ordinary sense. Even more, the sequence $l$ has no $r$-limits for any $r$ since it is unbounded from above (see Theorem 2.3 from [8]).

On the other hand, we see that the sequence $x$ is not statistically convergent because it does not satisfy the Cauchy condition for statistical convergence [21].

At the same time, $0 = stat\text{-}1\text{-}\lim l$ since $d(K) = 0$ where $K = \{n^2 \text{ for all } n = 1,2,3,\dots\}$.

**Lemma 3.2.** Statistical 0-convergence coincides with the concept of statistical convergence.

This result shows that statistical fuzzy convergence is a natural extension of statistical convergence.

**Lemma 3.3.** If $a = stat\text{-}\lim l$, then $a = stat\text{-}r\text{-}\lim l$ for any $r \geq 0$.

**Lemma 3.4.** If $a = stat\text{-}r\text{-}\lim l$, then $a = stat\text{-}q\text{-}\lim l$ for any $q > r$.

**Lemma 3.5.** If $a = stat\text{-}r\text{-}\lim l$ and $|b - a| = p$, then $b = stat\text{-}q\text{-}\lim l$ where $q = p + r$.

It is known that a subsequence of a fuzzy convergent sequence is fuzzy convergent [8]. However, for statistical convergence this is not true. Indeed, the sequence $h = \{i \; ; \; i = 1,2,3,\dots\}$ is a subsequence of the statistically fuzzy convergent sequence $l$ from Example 3.1. At the same time, $h$ is statistically fuzzy divergent.

However, if we consider dense subsequences of statistically fuzzy convergent sequences, it is possible to prove the following result.

**Theorem 3.1.** A sequence is statistically $r$-convergent if and only if any its statistically dense subsequence is statistically $r$-convergent.

*Proof*. *Necessity.* Let us take a statistically $r$-convergent sequence $l = \{a_i \; ; \; i = 1,2,3,\dots\}$ and a statistically dense subsequence $h = \{b_k \; ; \; k = 1,2,3,\dots\}$ of $l$. Let us also assume that $h$ statistically $r$-diverges. Then for any real number $a$, there is some $\varepsilon > 0$ such that $d(H_\varepsilon(a)) > 0$ where $H_\varepsilon(a) = \{k; |b_k - a| > r + \varepsilon \}$. As $h$ is a subsequence of $l$, we have $L_\varepsilon(a) \supseteq H_\varepsilon(a)$ where $L_\varepsilon(a) = \{i; |a_i - a| > r + \varepsilon\}$. Consequently, $d(H_\varepsilon(a)) > 0$ as the subsequence $h$ is statistically dense in $l$. Thus, $l$ is also statistically $r$-divergent.

*Sufficiency* follows from the fact that $l$ is a statistically dense subsequence of itself.

Theorem is proved.

A statistically $r$-convergent sequence contains not only dense statistically $r$-convergent subsequences, but also dense $r$-convergent subsequences.

**Theorem 3.2.** $a = stat\text{-}r\text{-}\lim l$ if and only if there exists an increasing index sequence $K = \{k_n \; ; \; k_n \in N, n = 1,2,3,\dots\}$ of the natural numbers such that $d(K) = 1$ and $a = r\text{-}\lim l_K$ where $l_K = \{a_i \; ; \; i \in K\}$.

*Proof. Necessity.* Suppose that $a = stat\text{-}r\text{-}\lim l$. Let us consider sets $L^{r,j}(a) = \{i \in N; |a_i - a| < r + (1/j) \}$ for all $j = 1,2,3,\dots$ By the definition, we have

(3.2)                          $L^{r,j+1}(a) \subseteq L^{r,j}(a)$

and as $a = stat\text{-}r\text{-}\lim l$, by Lemma 3.1, we have

(3.3)                          $d(L^{r,j}(a)) = 1$

for all $j = 1,2,3,\ldots$ Let us take some number $i_1$ from the set $L^{r,1}(a)$. Then, by (3.2) and (3.3), there is a number $i_2$ from the set $L^{r,2}(a)$ such that $i_1 < i_2$ and

$$(1/n)\,|\{i \in N; i \leq n ; |a_i - a| < r + 1/2\}| > 1/2 \qquad \text{for all } n \geq i_2 \,.$$

In a similar way, we can find a number $i_3$ from the set $L^{r,3}(a)$ such that $i_2 < i_3$ and

$$(1/n)\,|\{i \in N; i \leq n ; |a_i - a| < r + 1/3\}| > 2/3 \qquad \text{for all } n \geq i_3 \,.$$

We continue this process and construct an increasing sequence $\{i_j \in N, j = 1,2,3,\ldots\}$ of the natural numbers such that each number $i_j$ belongs to $L^{r,j}(a)$ and

(3.4)          $(1/n)\,|\{i \in N; i \leq n ; |a_i - a| < r + 1/j\}| > (j-1)/j \qquad \text{for all } n \geq i_j \,.$

Now we construct the increasing sequence of indices $K$ as follows:

(3.5)          $K = \{i \in N; 1 \leq i \leq i_1\} \cup (\cup_{j \in N} \{ i \in L^{r,j}(a); i_j \leq i \leq i_{j+1}\})$.

Then from (3.2), (3.4) and (3.5), we conclude that for all $n$ from the interval $i_j \leq n \leq i_{j+1}$ and all $j = 1,2,3,\ldots$, we have

(3.6)          $(1/n)\{k \in K; k \leq n\} = (1/n)\,|\{i \in N; i \leq n ; |a_i - a| < r + 1/j\}| > (j-1)/j.$

Hence it follows from (3.6) that $d(K) = 1$. Now let us denote $l_K = \{a_i ; i \in K\}$, take some $\varepsilon > 0$ and choose a number $j \in N$ such that $1/j < \varepsilon$. If $n \in K$ and $n \geq i_j$, then, by the definition of $K$, there exists a number $m \geq j$ such that $i_m \leq n \leq i_{m+1}$ and thus, $n \in L_{r,m}(a)$. Hence, we have

$$|a_n - a| < r + 1/j < r + \varepsilon$$

As this is true for all $n \in K$, we see that $a = stat\text{-}r\text{-}\lim l_K$

Thus, the proof of necessity is completed.

*Sufficiency.* Suppose that there exists an increasing index sequence $K = \{k_n ; k_n \in N; n = 1,2,3,\ldots\}$ of the natural numbers such that $d(K) = 1$ and $a = r\text{-}\lim l_K$ where $l_K = \{a_i ; i$ $\in K$}. Then there is a number $n$ such that for each $i$ from $K$ such that $i \geq n$, the inequality $|a_i - a| < r + \varepsilon$ holds. Let us consider the set

$$L_{r,\varepsilon}(a) = \{i \in N; |a_i - a| \geq r + \varepsilon\}$$

Then we have

$$L_{r,\varepsilon}(a) \subseteq N \setminus \{k_i; k_i \in N; i = n, n+1, n+2, \dots\}$$

Since $d(K) = 1$, we get $d(N \setminus \{k_i; k_i \in N; i = n, n+1, n+2, \dots\}) = 0$, which yields $d(L_{r,\varepsilon}(a)) = 0$ for every $\varepsilon > 0$. Therefore, we conclude that $a = stat\text{-}r\text{-}\lim l$.

Theorem is proved.

**Corollary 3.1** [32]. $a = stat\text{-}\lim l$ if and only if there exists an increasing index sequence $K = \{k_n; k_n \in N; n = 1,2,3,\dots\}$ of the natural numbers such that $d(K) = 1$ and $a = \lim l_K$ where $l_K = \{a_i; i \in K\}$.

**Corollary 3.2.** $a = stat\text{-}r\text{-}\lim l$ if and only if there exists a sequence $h = \{b_i; i = 1,2,3,\dots\}$ such that $d(\{i; a_i = b_i\}) = 1$ and $a = r\text{-}\lim h$.

**Corollary 3.3.** The following statements are equivalent:

(*i*) $a = stat\text{-}r\text{-}\lim l$.

(*ii*) There is a set $K \subseteq N$ such that $d(K) = 1$ and $a = r\text{-}\lim l_K$ where $l_K = \{a_i; i \in K\}$.

(*iii*) For every $\varepsilon > 0$, there exist a subset $K \subseteq N$ and a number $m \in K$ such that $d(K) = 1$ and $|a_n - a| < r + \varepsilon$ for all $n \in K$ and $n \geq m$.

We denote the set of all statistical $r$-limits of a sequence $l$ by $L_{r\text{-}stat}(l)$, that is,

$$L_{r\text{-}stat}(l) = \{a \in R; a = stat\text{-}r\text{-}\lim l\}.$$

Then we have the following result.

**Theorem 3.3.** For any sequence $l$ and number $r \geq 0$, $L_{r\text{-}stat}(l)$ is a convex subset of the real numbers.

*Proof*. Let $c, d \in L_{r\text{-}stat}$ ($l$), $c < d$ and $a \in [c, d]$. Then it is enough to prove that $a \in L_{r\text{-}stat}$ ($l$). Since $a \in [c, d]$, there is a number $\lambda \in [0, 1]$ such that. $a = \lambda c - (1-\lambda)d$. As $c, d \in L_{r\text{-}stat}$ ($l$), then for every $\varepsilon > 0$, there exist index sets $K_1$ and $K_2$ with $d(K_1) = d(K_2) = 1$ and the numbers $n_1$ and $n_2$ such that $|a_i - c| < r + \varepsilon$ for all $i$ from $K_1$ and $i \geq n_1$ and $|a_i - d| < r + \varepsilon$ for all $i$ from $K_2$ and $i \geq n_2$. Let us put $K = K_1 \cap K_2$ and $n = \max\{n_1, n_2\}$. Then $d(K) = 1$ and for all $i \geq n_1$ with $i$ from $K$, we have

$$
\begin{aligned}
| a_i - a | &= | a_i - \lambda c - (1-\lambda)d| \\
&= |(\lambda a_i - \lambda c) + ((1-\lambda)\, a_i - (1-\lambda)d)| \\
&\leq |\lambda a_i - \lambda c| + |(1-\lambda)\, a_i - (1-\lambda)d \,| \\
&\leq \lambda(r + \varepsilon) + (1-\lambda)(r + \varepsilon) = r + \varepsilon.
\end{aligned}
$$

So, by Theorem 3.1, we conclude $a = stat\text{-}r\text{-}\lim l$, which implies $a \in L_{r\text{-}stat}$ ($l$). Theorem is proved.

Lemmas 3.4 and 3.5 imply the following result.

**Proposition 3.1.** If $q > r$, then $L_{r\text{-}stat}$ ($l$) $\subset L_{q\text{-}stat}$ ($l$).

**Definition 3.2.** The quantity

$$\inf\{\, r \,;\, a = stat\text{-}r\text{-}\lim l \,\}$$

is called the upper statistical defect $\delta(a = stat\text{-}\lim l)$ of statistical convergence of $l$ to the number $a$.

**Proposition 3.2.** If $q = \inf\{r; a = stat\text{-}r\text{-}\lim l\}$, then $a = stat\text{-}q\text{-}\lim l$.

**Definition 3.3.** The quantity

$$\frac{1}{1 + \delta(a = stat\text{-}\lim l)}$$

is called the upper statistical measure $\mu(a = stat\text{-lim } l)$ of statistical convergence of $l$ to a number $a$.

The upper statistical measure of statistical convergence of $l$ defines the fuzzy set $\boldsymbol{L}_{stat}$ ($l$) = [$\boldsymbol{R}$, $\mu(a = stat\text{-lim } l)$], which is called the *complete statistical fuzzy limit* of the sequence $l$.

**Example 3.2.** We find the complete statistical fuzzy limit $\boldsymbol{L}_{stat}$ ($l$) of the sequence $l$ from Example 3.1. For this sequence and a real number $a$, we have

$$\mu(a = stat\text{-lim } l) = 1/(2 + |a|)$$

This fuzzy set $\boldsymbol{L}_{stat}$ ($l$) is presented in Figure 1.

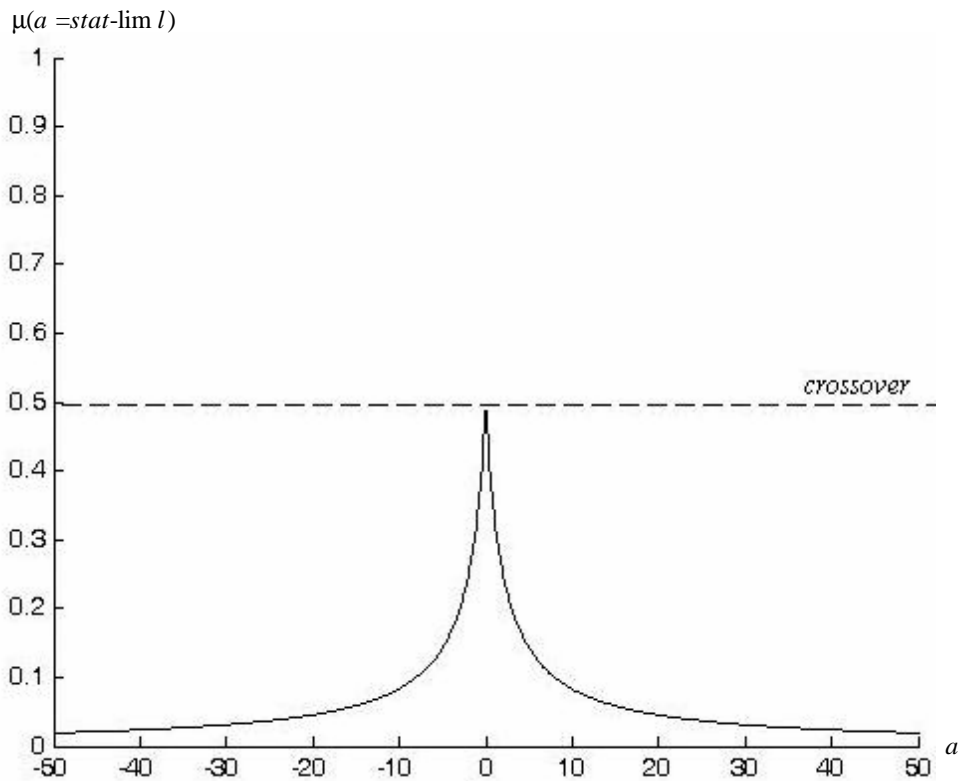

**Figure 1.** The complete statistical fuzzy limit of the sequence $l$ from Example 3.1.

**Definition 3.4** [36]**.** A fuzzy set [A, μ] is called *convex* if its membership function μ($x$) satisfies the following condition: μ($?x$ + (1 - $?$)$z$) = min{μ($x$), μ($z$)} for any $x$, $z$ and any number $? > 0$.

Then we have the following result.

**Theorem 3.4.** The complete statistical fuzzy limit $L_{stat}$ ($l$) = {$a$, μ($a$ = *stat*-$r$-lim $l$); $a$ ∈ $\boldsymbol{R}$} of a sequence $l$ is a convex fuzzy set.

*Proof*. Let $c$, $d$ ∈ $L_{r\text{-}stat}$ ($l$), $c < d$ and $a$ ∈ [$c$, $d$]. Then it is enough to prove that μ($a$ = *stat*-lim $l$) = μ( ($?c$ + (1 - $?$)$d$) = *stat*-lim $l$) = min{μ($c$ = *stat*-lim $l$), μ($d$ = *stat*-lim $l$ )}. This is equivalent to the inequality δ($a$ = *stat*-lim $l$) = δ( ($?c$ + (1 - $?$)$d$) = *stat*-lim $l$) ≤ max{δ($c$ = *stat*-lim $l$), δ($d$ = *stat*-lim $l$ )}.

Let us assume, for convenience, that $q$ =δ($c$ = *stat*-lim $l$) = $r$ = δ($d$ = *stat*-lim $l$ )}. Then by Lemma 3.4, $d$ = *stat*-$q$-lim $l$. Then by Theorem 3.3, $d$ = *stat*-$q$-lim $l$ as the set $L_{r\text{-}stat}$ ($l$) is convex. Thus, δ($a$ = *stat*-lim $l$) ≤ $q$ = max{δ($c$ = *stat*-lim $l$), δ($d$ = *stat*-lim $l$ )}.

Theorem is proved.

**Definition 3.5** [36]**.** A fuzzy set [A, μ] is called *normal* if sup μ($x$) = 1.

Theorem 3.4 allows us to prove the following result.

**Theorem 3.5.** The complete statistical fuzzy limit $L_{stat}$ ($l$) = {$a$, μ($a$ = *stat*-$r$-lim $l$); $a$ ∈ $\boldsymbol{R}$} of a sequence $l$ is a normal fuzzy set if and only if $l$ statistically converges.

Let $l$ = {$a_i$∈$\boldsymbol{R}$; $i$ = 1,2,3, …}and $h$ = {$b_i$∈$\boldsymbol{R}$; $i$ = 1,2,3, …}. Then their sum $l$ + $h$ is equal to the sequence {$a_i$ + $b_i$; $i$ = 1,2,3, …} and their difference $l$ - $h$ is equal to the sequence {$a_i$ - $b_i$; $i$ = 1,2,3, …}. Lemma 2.1 allows us to prove the following result.

**Theorem 3.6.** If $a$ = *stat*-$r$-lim $l$ and $b$ = *stat*-$q$-lim $h$, then:

(a) $a + b = stat\text{-}(r+q)\text{-}\lim(l+h)$ ;

(b) $a - b = stat\text{-}(r+q)\text{-}\lim(l - h)$ ;

(c) $ka = stat\text{-}(\,|k|\cdot r)\text{-}\lim\,(kl)$  for any $k \in \boldsymbol{R}$  where $kl = \{ka_i \,;\, i = 1,2,3,\dots\}$.

**Corollary 3.3** [32]. If $b = stat\text{-}\lim l$ and  $c = stat\text{-}\lim h$, then:

(a) $a + b = stat\text{-}\lim(l+h)$ ;

(b) $a - b = stat\text{-}\lim(l - h)$ ;

(c) $ka = stat\text{-}\lim\,(kl)$  for any $k \in \boldsymbol{R}$.

An important property in calculus is the Cauchy criterion of convergence, while an important property  in neoclassical analysis is the extended Cauchy criterion of fuzzy convergence. Here we find an extended statistical Cauchy criterion for statistical fuzzy convergence.

**Definition 3.6.**  A sequence $l$ is called *statistically r-fundamental* if for any  $\varepsilon > 0$ there is $n \in \boldsymbol{N}$ such that $d(L_{n,r,\varepsilon}) = 0$ where $L_{n,r,\varepsilon} = \{\, i \in \boldsymbol{N};\, i \leq n$ and  $|a_i - a_n| \geq r + \varepsilon \,\}$.

**Definition 3.7.**  A sequence $l$ is called *statistically fuzzy fundamental* if it is statistically $r$-fundamental for some $r \geq 0$.

**Lemma 3.6.**  If  $r \leq p$,  then any statistically $r$-fundamental sequence is statistically $p$-fundamental.

**Lemma 3.7.**  A sequence $l$ is a statistically Cauchy sequence if and only if it is statistically 0-fundamental.

This result shows that the property to be a statistically fuzzy fundamental sequence is a natural extension of the property to be a statistically Cauchy sequence.

Using the similar technique as in proof of Theorem 3.2, one can obtain the following result.

**Theorem 3.7.** A sequence $l$ is statistically $r$-fundamental if and only if there exists an increasing index sequence $K = \{k_n\,;\, k_n \in \mathbf{N}\,;\, n = 1,2,3,\dots\}$ of the natural numbers such that $d(K) = 1$ and the subsequence $l_K$ is $r$-fundamental, that is, for every $\varepsilon > 0$ there is a number $i$ such that $|a_{kn} - a_{ki}| < r + \varepsilon$ for all $n \geq i$.

**Corollary 3.4.** A sequence $l$ is statistically fuzzy fundamental if and only if there exists a statistically dense subsequence $u$ such that $u$ is fuzzy fundamental.

Theorem 3.6 yields the following result.

**Corollary 3.5** [21]. A sequence $l$ is a statistically Cauchy sequence if and only if there exists an increasing index sequence $K = \{k_n\,;\, k_n \in \mathbf{N}\,;\, n = 1,2,3,\dots\}$ of the natural numbers such that $d(K) = 1$ and the subsequence $l_K$ is a Cauchy sequence.

**Theorem 3.8** *(The Extended Statistical Cauchy Criterion).* A sequence $l$ has a statistical $r$-limit if and only if it is statistically $r$-fundamental.

*Proof. Necessity.* Let $a = stat$-$r$-lim $l$. Then by the definition, for any $\varepsilon > 0$, we have $d(L_{r,\varepsilon}(a)) = 0$, in other words, $\lim_{n\to\infty} (1/n)\,|\{i \in \mathbf{N}\,;\, i \leq n\,;\, |a_i - a| \geq r + \varepsilon/2\}| = 0$. This implies that given $\varepsilon > 0$, we find $n \in \mathbf{N}$ such that for any $i > n$, we have $|a_i - a_n| \leq |a - a_i| + |a - a_n|$. Consequently, $d(L_{n,r,\varepsilon}) \leq d(L_{r,\varepsilon/2}(a)) + d(L_{r,\varepsilon/2}(a)) = 0$, i.e., $d(L_{n,r,\varepsilon}) = 0$ Thus, $l$ is a statistically $r$-fundamental sequence.

*Sufficiency.* Assume now that $l$ is a statistically $r$-fundamental sequence. Then, by Theorem 3.7, we conclude that there is an $r$-convergent $u = \{u_i;\, i = 1,2,3,\dots\}$ such that $d(\{i\,;\, a_i = u_i\}) = 1$. We denote the $r$-limit of $u$ is $a$. Now let $A = \{i \in \mathbf{N};\, i \leq n\,;\, a_i \neq u_i\}$ and $B = \{i \in \mathbf{N}\,;\, i \leq n\,;\, |u_i - a| \geq r + \varepsilon\}$. Then observe that $d(A) = d(B) = 0$. On the other hand, since for each $n$

$$L_{r,\varepsilon}(a) = \{i \in N; i \leq n; \mid a_i - a \mid \geq r + \varepsilon\} \subseteq A \cup B,$$

we have $d(L_{r,\varepsilon}) = 0$, which gives *stat*-*r*-lim $l = a$.

The proof is completed.

Theorem 3.8 directly implies the following results.

**Corollary 3.6** *(The General Fuzzy Convergence Criterion).* The sequence *l* statistically fuzzy converges if and only if it is statistically fuzzy fundamental.

**Corollary 3.7** (*The Statistical Cauchy Criterion*) [21]. A sequence *l* statistically converges if and only if it is statistically fundamental, i.e., for any $\varepsilon > 0$ there is $n \in N$ such that $d(L_{n,r,\varepsilon}) = 0$ where $L_{n,r,\varepsilon} = \{i \in N; \mid a_i - a_n \mid \geq \varepsilon\}$.

**Corollary 3.8** (*The Cauchy Criterion*). The sequence *l* converges if and only if it is fundamental.

## 4. Fuzzy convergence in statistics and statistical fuzzy convergence

In Section 2, we found relations between statistical convergence and convergence of statistical characteristics (such as mean and standard deviation). However, when data are obtained in experiments, they come from measurement and computation. As a result, we never have and never will be able to have absolutely precise convergence of statistical characteristics. It means that instead of ideal classical convergence, which exists only in pure mathematics, we have to deal with fuzzy convergence, which is closer to real life and gives more realistic models. That is why in this section, we consider relations between statistical fuzzy convergence and fuzzy convergence of statistical characteristics.

Let $l = \{a_i ; i = 1,2,3,\ldots\}$ be a bounded sequence, i.e., there is a number $m$ such that $|a_i| < m$ for all $i \in \mathbf{N}$. This condition is usually true for all sequences generated by measurements or computations.

**Theorem 4.1.** If $a = stat\text{-}r\text{-}\lim l$, then $a = r\text{-}\lim \mu(l)$ where $\mu(l) = \{\mu_n = (1/n) \sum_{i=1}^{n} a_i;$ $n = 1, 2, 3, \ldots\}$.

*Proof.* Since $a = stat\text{-}r\text{-}\lim l$, for every $\varepsilon > 0$, we have

(4.1)                  $\lim_{n \to \infty} (1/n) |\{ i \leq n, i \in \mathbf{N}; | a_i - a| \geq r + \varepsilon \}| = 0.$

If $|a_i| < m$ for all $i \in \mathbf{N}$, then there is a number $k$ such that $|a_i - a| < k$ for all $i \in \mathbf{N}$. Namely, $|a_i - a| \leq |a_i| + |a| \leq m + |a| = k$. Taking the set $L_{n,r,\varepsilon}(a) = \{i \in \mathbf{N}; i \leq n$ and $| a_i - a | \geq r + \varepsilon\}$, denoting $|L_{n,r,\varepsilon}(a)|$ by $u_n$, and using the hypothesis $|a_i| < m$ for all $i \in \mathbf{N}$, we have the following system of inequalities

$$|\mu_n - a| = |(1/n) \sum_{i=1}^{n} a_i - a|$$

$$\leq (1/n) \Sigma_{i=1}^{n} |a_i - a|$$

$$\leq (1/n) (ku_n + (n - u_n)(r + \varepsilon))$$

$$\leq (1/n) (ku_n + n(r + \varepsilon))$$

$$= r + \varepsilon + (1/n) (ku_n).$$

From the equality (4.1), we get, for sufficiently large $n$, that the inequality $|\mu_n - a| \leq r + 2\varepsilon$ holds because the number $k$ is fixed and the sequence $\{(1/n)u_n ; n = 1,2,3,\ldots\}$ converges to zero. Thus, $a = r\text{-}\lim \mu(l)$.

Theorem is proved.

**Remark 4.1.** Statistical *r*-convergence of a sequence does not imply *r*-convergence of this sequence even if all elements are bounded as the following example demonstrates.

**Example 4.1.** Let us consider the sequence $l = \{a_i ; i = 1,2,3,\ldots\}$ elements of which are

$$a_i = \begin{cases} (-1)^i \cdot 1000 & \text{when } i = n^2 \text{ for all } n = 1,2,3,\ldots \\ \\ (-1)^i & \text{otherwise.} \end{cases}$$

By the definition, $0 = stat\text{-}1\text{-}\lim l$ since $d(K) = 0$ where $K = \{n^2 \text{ for all } n = 1,2,3,\ldots\}$. At the same time, this sequence does not have 1-limits.

**Corollary 4.1.** If sequence $l$ is statistically fuzzy fundamental, then the sequence of its partial means is fuzzy fundamental.

Finally, we get the following result.

**Theorem 4.2.** If $a = stat\text{-}r\text{-}\lim l$ and there is a number $m$ such that $|a_i| < m$ for all $i = 1,2,\ldots$, then $0 = [2pr]^{\frac{1}{2}}\text{-}\lim \sigma(l)$ where $p = \max \{m^2 + |a|^2, m + |a|\}$.

*Proof*. We will first show that $\lim \sigma^2(l) = 0$. By the definition, $\sigma_n^2 = (1/n) \Sigma_{i=1}^{n} (a_i - \mu_n)^2 = (1/n) \Sigma_{i=1}^{n} (a_i)^2 - \mu_n^2$. Thus, $\lim \sigma^2(l) = \lim_{n\to\infty} (1/n) \Sigma_{i=1}^{n} (a_i)^2 - \lim_{n\to\infty} \mu_n^2$. Since $|a_i| < m$ for all $i \in N$, there is a number $p$ such that $|a_i^2 - a^2| < p$ for all $i \in N$. Namely, $|a_i^2 - a^2| \le |a_i|^2 + |a|^2 < m^2 + |a|^2 < \max \{m^2 + |a|^2, m + |a|\} = p$. Taking the set $L_{n,r,\varepsilon}(a) = \{i \in N; i \le n \text{ and } |a_i - a| \ge r + \varepsilon\}$, denoting $|L_{n,r,\varepsilon}(a)|$ by $u_n$, and using the hypothesis $|a_i| < m$ for all $i \in N$, we have the following system of inequalities:

$$|\sigma^2_n| = |(1/n) \Sigma_{i=1}^{n} (a_i)^2 - \mu_n^2|$$

$$= |(1/n) \Sigma_{i=1}^{n} (a_i^2 - a^2) - (\mu_n^2 - a^2)|$$

$$\le (1/n) \Sigma_{i=1}^{n} |a_i^2 - a^2| + |\mu_n^2 - a^2|$$

$$< (p/n) \Sigma_{i=1}^{n} |a_i - a| + |\mu_n - a| \, |\mu_n + a|$$

$$< (p/n) (u_n + (n - u_n)(r + \varepsilon)) + |\mu_n - a| (|\mu_n| + |a|)$$

$$\le (p/n) (u_n + n (r + \varepsilon)) + |\mu_n - a| ((1/n)\Sigma_{i=1}^{n} |a_i| + |a|)$$

$$< p \ (u_n \ /n) + p \ (r + \varepsilon) + p \ |\mu_n - a|.$$

Now by hypothesis and Theorem 4.1, we have $a = r\text{-}\lim \mu(l)$. Also, by (4.1), we have $\lim (u_n \ /n) = 0$. Then, for every $\varepsilon > 0$ and sufficiently large $n$, we have

(4.2)                              $|\sigma^2{}_n \ | < p \ \varepsilon + p \ (r + \varepsilon) + p \ (r + \varepsilon) = 2pr + 3p\varepsilon.$

As $(x + y)^{\frac{1}{2}} \le x^{\frac{1}{2}} + y^{\frac{1}{2}}$ for any $x, y > 0$, it follows from (4.2) that

$$|\sigma_n \ | \le [2pr]^{\frac{1}{2}} + (3p\varepsilon)^{\frac{1}{2}},$$

which yields that $0 = [2pr]^{\frac{1}{2}}\text{-}\lim \sigma(l)$.

The proof is completed.

Fuzzy convergence of partial standard deviations implies corresponding fuzzy convergence of partial variances.

**Corollary 4.2.** If $a = stat\text{-}r\text{-}\lim l$ and there is a number $m$ such that $| \ a_i| < m$ for all $i = 1, 2, 3, \ldots$, then $0 = [2pr]\text{-}\lim \sigma^2(l)$ where $p = \max \ \{m^2 + |a|^2, m + |a|\}$.

## 5. Conclusion

We have developed the concept of statistical fuzzy convergence and studied its properties. Relations between statistical fuzzy convergence and fuzzy convergence are considered in Theorems 3.1 and 3.2. Algebraic structures of statistical fuzzy limits are described in Theorem 3.5. Topological structures of statistical fuzzy limits are described in Theorems 3.3 and 3.4. Relations between statistical convergence, ergodic systems, and convergence of statistical characteristics, such as the mean (average), and standard deviation, are studied in Sections 2 and 4. Introduced constructions and obtained results open new directions for further research.

It would be interesting to develop connections between statistical fuzzy convergence and fuzzy dynamical systems [9], introducing and studying ergodic systems.

It would also be interesting to study statistically fuzzy continuous functions, taking as the base the theory of fuzzy continuous functions [7] and the theory of statistically continuous functions [10].

The theory of regular summability method is an important topic in functional analysis (see, for instance, [6, 22]). In recent years it has been demonstrated that statistical convergence can be viewed as a regular method of series summability. In particular, Connor showed [11] that statistical convergence is equivalent to the strong Cesaro summability in the space of all series with bounded elements. Similar problems are studied in [12]. Thus, it would also be interesting to study summability of series and relations between statistical fuzzy convergence and summability methods.